\documentclass[10pt,draftcls,onecolumn]{IEEEtran}

\usepackage{amsfonts}
\usepackage{amsmath}
\usepackage{amssymb}
\usepackage{algorithm,algorithmic}
\usepackage{booktabs}
\usepackage{graphicx}
\usepackage{xcolor}
\usepackage{enumerate}
\usepackage{amsbsy}
\usepackage{booktabs}
\usepackage{subcaption}
\usepackage{url}

\newtheorem{theorem}{Theorem}

\newtheorem{remark}{Remark}

\newtheorem{definition}{Definition}

\newcommand{\real}{\text{Re}}
\newcommand{\imag}{\text{Im}}
\renewcommand{\i}{{\text{i}}}

\newcommand{\ped}[1]{{_{\mathrm{#1}}}}

\newcommand{\uu}{\mathbf{u}}
\newcommand{\xx}{\mathbf{x}}
\newcommand{\VV}{\mathbf{V}}
\newcommand{\WW}{\mathbf{W}}
\newcommand{\dd}{\boldsymbol{\delta} }
\renewcommand{\aa}{\boldsymbol{\alpha}}
\newcommand{\DD}{\boldsymbol{\Delta} }

\definecolor{green}{rgb}{0,0.7,0}

\def\tSOC{\theta\ped{SwC}}

\def\V{\mbox{\rm Viol}}

\newcommand{\Real}[1]{\mathbb{R}^{#1}}

\IEEEoverridecommandlockouts

\newcommand{\FD}[1]{{#1}}

\title{\LARGE \bf{Probabilistically Robust AC Optimal Power Flow}}
\begin{document}

\author{Mohammadreza~Chamanbaz,~\IEEEmembership{Member,~IEEE,}
Fabrizio Dabbene,~\IEEEmembership{Senior~Member,~IEEE,}, and
Constantino Lagoa,~\IEEEmembership{Member,~IEEE}%
\thanks{M. Chamanbaz is with iTrust Center for Research in Cyber Security, Singapore University of Technology and Design, 8 Somapah Road Singapore 487372, (E-mail: Chamanbaz@sutd.edu.sg, Chamanbaz@u.nus.edu).}
\thanks{F. Dabbene is with CNR-IEIIT, Corso Duca degli Abruzzi 24, 12129 Torino, Italy (E-mail:fabrizio.dabbene@ieiit.cnr.it).  }
\thanks{C. Lagoa is with the Department of Electrical Engineering and Computer Science, Penn State University, University Park, PA 16802 USA, (E-mail:  lagoa@engr.psu.edu).}
\thanks{This work was supported in part by the CNR International Joint Lab COOPS, the National Science Foundation under grant \mbox{CNS-1329422}, and the Singapore National Research Foundation (NRF) grant
under the ASPIRE project, grant No NCR-NCR001-040.}}

\markboth{IEEE Transactions on Control of Network Systems}{}%
	
\maketitle
	
\begin{abstract}
The increasing penetration of renewable energy resources, paired with the fact that load can vary significantly, introduce a high degree of uncertainty in the behavior of modern power grids. Given that classical dispatch solutions are ``rigid," their performance in such an uncertain environment is in general far from optimal. For this reason, in this paper, we consider AC optimal power flow (AC-OPF) problems in the presence of uncertain loads and (uncertain) renewable energy generators. The goal of AC-OPF design is to guarantee that controllable generation is dispatched at minimum cost, while satisfying constraints on generation and transmission for almost all realizations of the uncertainty.
 We propose an approach based on a randomized technique recently developed, named \textit{scenario with certificates}, 
 which allows us to tackle the problem without the conservative parameterizations on the uncertainty used in currently available approaches.
 The proposed solution can exploit the usually available probabilistic description of the uncertainty and variability, and provides solutions with \textit{a-priori} probabilistic guarantees on the risk of violating the constraints on generation and transmission.
\end{abstract}

\section{Introduction}
Modern power grids are characterized by increasing penetration of renewable energy sources, such as solar photovoltaic and wind power. 
This trend is expected to increase in the near future, as also testified by strict commitments to large renewable power penetration being made by major countries worldwide; e.g., see  \cite{SET2016,GWEC2016}. 
While the advantages of renewable energy in terms of environmental safeguard are indisputable, its introduction does not come without a cost. Indeed, renewable energy generation technologies are highly variable and not fully dispatchable, thus imposing novel challenges to the existing power system operational paradigm. As discussed in e.g. \cite{Bienstock2014}, when  uncontrollable resources fluctuate, classical optimal power flow (OPF) solutions can provide very inefficient power generation policies, that result in line overloads and, potentially,  cascading outrages.

Despite the increasingly larger investments, which are costly and subject to several regulatory and policy limitations, 
 power outages due to the uncertainty introduced by renewable power generation still occur frequently.
This situation clearly shows that a strategy based only on investments in technological improvements of the transmission lines and controllable generation capacity---as those discussed e.g.\ in \cite{
Conejo2010}---is not 
sufficient anymore. Instead, radically new dispatch philosophies need to be devised, able to cope with the increasing \textit{uncertainty}, due to unpredictable fluctuations in renewable output and time-varying loads. 

Indeed, classical OPF dispatch is typically computed based on simple predictions of expected loads and generation levels for the upcoming time window. Although these predictions can be fairly precise for the case of traditional generators and loads, they may be highly unreliable in the case of renewable generators, thus explaining its failure in these latter situations. 

It follows that one of the major problems in today's power grids is the following: \textit{Given the high level of uncertainty introduced  by renewable energy sources, design a dispatch policy that {\it i)} minimizes generation costs and {\it ii)} does not violate generation and transmission constraints for all admissible values of renewable power and variable demand. In other words, one would like to design an optimal dispatch policy that is robust against uncertainty.}

However, such robust policies might be very conservative. Power networks can tolerate temporary violations of their generation and transmission constraints. Moreover, being robust with respect to any possible value of the uncertainty may lead to very inefficient policies, since some scenarios are very unlikely. Hence, in this paper we take a different approach: instead of requiring that the network constraints are satisfied for all possible values of uncertainty, we allow for a small well-defined risk of constraint violation.
More precisely, we start by assuming that one can adapt the power generated by conventional generators, based on real-time information about current values of renewable power generation and demand. Under this assumption, we aim at minimizing cost of power generation while meeting demand and allowing a small well defined risk of violating network constraints. 
We re-formulate the problem by highlighting the fundamental difference between control and state variables. Then, to approximate the solution of the resulting (complex) optimization problem we leverage recent results on convex relaxations of the optimal power flow problem~\cite{Lavaei:2012}, 
 and show that such a relaxation is exact for a restrictive class of networks even in the presence of uncertainty. 
Finally, together with the above mentioned relaxation, we  use a novel way of addressing probabilistic robust optimization problems known as \emph{scenario with certificates} (SwC) \cite{SWC-TAC} to optimize generation cost subject to a small risk of network failure. In this way, we derive a convex problem that i) is efficiently solvable, ii) provides a good approximation of the optimal power generation under the above mentioned risk constraints and iii) for some special classes of networks, it provides an exact solution.

\vspace{-0.1 in}
\subsection{Previous Results} 

\subsubsection{Nominal OPF}
The optimal power flow problem is known to be \mbox{NP-hard} even in the absence of uncertainty; see e.g.,~\cite{NP-Hard-1,FERC-1,NP-Hard-2} and references therein. Hence, several numerical approaches propose approximations of OPF, based e.g.\ on Newton methods~\cite{DommelTinney1968},  interior point based methods~\cite{YanQuintana1999} or global optimizations heuristics \cite{Abido2002,LaiMaYokoyamaEtAl1997}. 
In particular, several approximations have been introduced to recover convexity and make the problem numerically tractable. The most common is the so-called DC approximation, see \cite{Coffrin:2014} and references therein, in which the AC OPF problem is linearized. However, the solution is in general sub-optimal and, more importantly,  it may not be feasible, in the sense that it may not satisfy the original nonlinear power flow equations. 
Also, as noted in \cite{Coffrin:2014},  the fact that DC approximation fixes voltage magnitudes and ignores reactive power, makes the solution not applicable in several important practical situations.

Motivated by the above considerations, semidefinite programming (SDP) relaxations of the AC OPF problem have been recently introduced to alleviate the computational burden, in \cite{Jabr:06} for radial networks and in \cite{Bai:08} for general networks---see \cite{Low:2014a} for a detailed review. These convex relaxations 
have recently received renewed interest thanks to \cite{Lavaei:2012},  that analyzes the optimality properties of the relaxation, showing how in many practical situation these relaxations turn out to be non-conservative---in the sense that the solution of the relaxed problem coincides with that of the original non-convex one. 
This result has sparked an {interesting literature} analyzing specific cases in which the SDP relaxations can be proven to be exact, see e.g. \cite{Farivar:13}, and showing how graph sparsity can be exploited to simplify the SDP relaxation of OPF, see e.g.\ \cite{Bai:11}.

However, while these works have reached a good level of maturity, (see for instance the recent two-part tutorial \cite{Low:2014a,Low:2014b}), 
there is still very little research analyzing if and how these relaxations could be extended to the problem considered in this paper, namely AC OPF in the presence of possibly large and unpredictable uncertainties, as detailed next. 

\subsubsection{DC-based approaches to uncertain OPF}

Most literature on uncertain OPF gets around the nonlinearities by recurring  to  the DC-based approximation. 
 These assumptions reduce the optimization problem to a quadratic program subject to uncertain linear equalities and inequalities, which still represents a challenging problem, at least for general probability distributions. 
First approaches in this direction have been based on scenario-tree generation methods, see for instance \cite{Yong2000}. These techniques suffer from severe computational complexity limitations, and do not offer theoretical guarantees on the probability of satisfaction of the constraints of the found solution.

In the case uncertainties are assumed  to be Gaussian, the problem can be written in closed form \cite{Andersson2013}, or is amenable to a second-order cone-program \cite{Bienstock2014},  for which efficient solutions exist. This approach has been extended to ambiguous densities introducing the so-called robustified chance constraints in \cite{Lubin2016}. Similar ideas are at the basis of the approaches in \cite{Morari2015,Xie2018},  which consider distributionally robust approaches. That is, the 
solution has to be valid for all uncertainties whose probability density  functions (pdfs) belong to a 
family of distribution functions sharing the same mean and variance, leading to the so-called conditional value at risk (CVar) optimization, which is again a convex problem.

The problem becomes much more difficult for general uncertainty distributions (indeed, it is known that the distribution of wind power is not Gaussian \cite{hodge_wind_2012}). In this case, a very promising approach is the application of recent results based on the so-called scenario approach, which are based on random  generation  of uncertainty samples. This is the approach followed for instance in \cite{ETH-PMAPS,ETH-TPS,ETH-Springer,ETH-PSCC}.

\subsubsection{Convex relaxation approaches for uncertain OPF}

The limitations inherent to the DC approximation have motivated a few recent approaches to uncertain OPF, which employ more sophisticated relaxations able to capture also the reactive components of the power equations.
The work \cite{Perninge2013}  makes use of second-order approximations of the stability boundary to approximate the probability of line violations. 
In \cite{ETH-AC} the SDP based relaxation introduced in \cite{Lavaei:2012} is exploited to derive a solution based on the scenario approach.
However, in order to guarantee solvability of the robust problem, the authors need to parameterize the dependence of the state variables on the uncertainty, see Remark \ref{crac-remark} for a detailed discussion.

Inspired by these recent works, in this paper we provide a less conservative approach, that does not require dependent variable parameterization. Instead, we exploit a recently developed approach to probabilistic robust optimization \emph{scenario with certificates} \cite{SWC-TAC} to develop less conservative, efficient relaxations of the OPF in the presence of uncertainty. As compared to previous results, in this paper we provide an approach to  the robust optimal power flow problem that: i) is exact for a (small) class of networks; 
ii) is computationally more efficient than other approaches in the literature and
iii) provides solutions leading to low risk of network failure even for networks that do not satisfy the theoretical requirements, as exemplified in the simulations provided.

\vspace{-0.1 in}
\subsection{The Sequel}
The paper is organized as follows: In Section~\ref{sec:prob formulation} we precisely formulate the AC-OPF problem. We also describe the adopted dispatch policy to be used in real time to cope with uncertainty. Then, we introduce a precise formulation of the Robust AC-OPF problem, which  divides the set of optimization variables into two distinct classes: i) \emph{independent/control variables} (those that can be controlled by the operator), and ii) \textit{dependent/state variables} (representing the state of the system). 
This reformulation represents one of the main contributions of our work, allowing to represent in a non-conservative way the problem of AC-OPF design in the presence of uncertainty. 
Also, in Section~\ref{sec:relaxation}, we show how this reformulation can be combined with the convex relaxation of~\cite{Lavaei:2012}, and we formally prove that this allows to recover its properties in terms of exactness of the relaxation in some special cases. 
Moreover, in Section \ref{sec:scenario} we show how this new formulation directly translates in the scenario with certificates paradigm. Numerical examples illustrating the performance of the proposed approach are provided in Section~\ref{sec:examples}. Finally, in Section~\ref{sec: conclusions}, some concluding remarks are presented.

\section{Optimal Power Flow Allocation Problem under Uncertainty}
\label{sec:prob formulation}

In this section, we briefly summarize the adopted power flow model, which  takes into explicit account load uncertainties and variable generators, and we formalize the ensuing optimization problem, arising from the necessity of robustly guaranteeing that safety limits are not exceeded. 

{ 

\vspace{-0.1 in}
\subsection{Robust AC-OPF}
We consider a power network with  graph representation 
$
\mathbb{G}=\{\mathcal{N},\mathcal{L}\},
$
where $\mathcal{N}\doteq\{1,2,\ldots,n\}$ denotes the set of buses (which can be represented as nodes of the graph), and $\mathcal{L}\subseteq\mathcal{N}\times\mathcal{N}$ denotes the set of electrical lines connecting the different buses in the network (represented as edges of the graph). The set of conventional generator buses is denoted by $\mathcal{G}\subseteq\mathcal{N}$, and its  cardinality is $n_g$. As a convention, it is assumed that the bus indices are ordered so that the first $n_g$ buses are generators, i.e. $\mathcal{G}\doteq\{1,2,\ldots,n_g\}$.
Each generator (we assume for ease of notation that no more than one generator is present on each generator bus) connected to the bus $k\in\mathcal{G}$ provides  complex power 
$
\bar{P}^G_{k}+ \bar{Q}^G_{k} \i$,
where $\bar{P}^G_{k}$ is the active power generated by  $k$-th generator, and $\bar{Q}^G_{k}$ is the corresponding reactive power. Goal of the network manager is to guarantee that the output of generators is such that the network operates safely and, if possible, at minimal cost.

We now elaborate on the different constraints that should be satisfied in order 
for the network to operate safely. 
In the framework considered in this paper, we assume both the existence of 
renewable energy sources and uncertain load. 
A renewable energy generator connected to  bus $k\in\mathcal{N}$ provides an 
uncertain complex power 
\begin{equation}\label{eq: renewable uncertainty}
P^R_{k}(\delta^R_{k})+Q^R_{k}(\delta^R_{k})\i=P^{R,0}_{k}+Q^{R,0}_{k}\i+\delta^R_{k},
\end{equation}
with $P^{R,0}_{k}+Q^{R,0}_{k}\i$ being the nominal (predicted) power generated 
by the \FD{RES}\footnote{We use the convention that $P^R_{k}=0$, 
$Q^R_{k}=0$ if no renewable generator is connected to node $k$.}, and 
$\delta^R_{k}\in\boldsymbol{\Delta}^R_{k}\subset\mathbb{C}$ representing an 
uncertain complex fluctuation, which mainly depends on the environmental 
conditions, such as wind speed in the case of wind generators. 

The uncertain demand in bus $k\in\mathcal{N}$ is also represented as
\begin{equation}\label{eq: load uncertainty }
P^L_{k}(\delta^L_{k})+Q^L_{k}(\delta^L_{k})\i=P^{L,0}_{k}+Q^{L,0}_{k}\i+\delta^L_{k}
\end{equation}
where $P^{L,0}_{k}$ and $Q^{L,0}_{k}$ denote the expected active and reactive 
load and  $\delta^L_{k}\in\boldsymbol{\Delta}^L_{k}\subset\mathbb{C}$ is the 
complex fluctuation in the demand at bus $k\in\mathcal{N}$.
The support set is the point $\{0\}$ if no uncertainty (i.e. no renewable 
generator or variable load) is present in bus $k$.
To simplify the notation, we collect the different sources of uncertainty by 
introducing  \textit{uncertainty vector} 
$
\dd \doteq
[\delta^L_{1}\,\cdots\,\delta^L_{n}\,\delta^R_{1}\,\cdots\,\delta^R_{n}]^T,
$
which varies in the set 
$
\DD \doteq \boldsymbol{\Delta}^L_{1} \times \cdots 
\times\boldsymbol{\Delta}^L_{n} \times
\boldsymbol{\Delta}^R_{1} \times \cdots \times \boldsymbol{\Delta}^R_{n}. 
$

To deal in a rigorous way with the uncertainties introduced above, we adopt a 
modification of the so-called frequency control (primary and secondary 
control), similar to that discussed in  \cite{Bienstock2014} in the 
context of DC-OPF. In classical OPF, 
this approach is used to distribute to the generators the difference between 
real-time (actual) and predicted demand,  through some coefficients which are 
generator specific. However, these coefficients are in general decided \textit{a-priori} 
in an ad-hoc fashion. 
This approach worked well in cases where the amount of power mismatch  was not 
significant; however, once renewable generators are in the power network, this 
difference may become large, thus leading to line overloads in the network.
The approach presented in  \cite{Bienstock2014} specifically incorporates 
these distribution parameters in the OPF optimization problem. In the setup 
proposed in this paper, we follow a similar strategy, and formally introduce a 
\textit{deployment vector}
$
\aa\doteq[\alpha_1,\ldots,\alpha_{n_g}]^T,
$
with 
$\sum_{k\in\mathcal{G}}\alpha_k=1$, $\alpha_k \geq 0$ for all $k \in 
\mathcal{G}$,
whose purpose is to distribute among the available generators 
the power mismatch created by the uncertain generators and loads. 
To the best of our knowledge, the use of a deployment vector was originally 
introduced in the context of DC-approximations in 
\cite{ETH-PMAPS} and \cite{ETH-Springer}, where it is referred to as ``distribution vector".
The same concept is present in many other works under different terminologies, such as  
 ``participation factor" in \cite{Jabr2013}, ``corrective control" in \cite{Jabr2015}, or ``affine control" in 
  \cite{Morari2013}. 

During operation, the active generation output of each generator is adapted 
according to the  realizations of the uncertain loads and renewable power (which are 
assumed to be measured on-line) as follows
\begin{align}\label{eq: affine control}
\bar{P}^G_{k} &= 
P^G_{k}+\alpha_k\left(\sum_{j\in\mathcal{N}}\real\{\delta^L_{j}\}-\sum_{k\in\mathcal{R}}\real\{\delta^R_{k}\}\right)\\
&= P^G_{k}+\alpha_k \mathbf{s}^T \real\{\dd\},\quad \forall k \in \mathcal{G} 
\nonumber
\end{align}
with $\mathbf{s}^T\doteq[\mathbf{1}_n^T,\,-\mathbf{1}_n^T]$.
Now, consider the line $(l,m)\in \mathcal{L}$, which is the line connecting buses $l$ and $m$.  Let  $y_{lm}$ be the (complex) admittance of the line and $V_k=|V_k|\angle{\theta_k}$ be the (complex) voltage at bus $k$, with magnitude $|V_k|$
and phase angle $\theta_k$. Then, the following \textit{Balance Equations} should be satisfied at all times
{\small \begin{subequations}\label{eq:unc_bal}
		\begin{align}
		\label{eq:unc_act_bal}
		&P^G_{k}+\alpha_k \mathbf{s}^T \real\{\dd\}  + P^R_{k}(\dd) -  
		P^L_{k}(\dd) = \\ \nonumber
		&\hspace{1 in}\sum_{l\in\mathcal{N}_k} 
		\real\left\{V_k(V_k-V_l)^*y_{kl}^*\right\}, \quad \forall k \in 
		\mathcal{N}\\
		\label{eq:unc_reac_bal}
		&Q^G_{k}  +  Q^R_{k}(\dd) - Q^L_{k}(\dd) = \\ \nonumber & \hspace{1 
		in}\sum_{l\in\mathcal{N}_k} 
		\imag\left\{V_k(V_k-V_l)^*y_{kl}^*\right\},\quad \forall k \in 
		\mathcal{N},\\
		&\mathbf{1}^T\aa=1 \label{eq:unc_bal alpha}
		\end{align}
\end{subequations}}
where $\mathcal{N}_k$ is the set of all neighboring buses directly connected to bus $k$.
Also, at generator bus $k\in \mathcal{G}$, one has the so-called \textit{Power Constraints}, which restrict the active and reactive power, and have the form
\begin{subequations}\label{eq:unc power constraints}
	\begin{align}
	\label{eq:unc power constraints active generation bound}
	&P_{k\,\min}\leq P^G_{k}+\alpha_k \mathbf{s}^T \real\{\dd\} \leq 
	P_{k\,\max},\quad \forall k \in \mathcal{G} \\
	\label{eq:unc power constraints reactive generation bound}
	&Q_{k\,\min}\leq Q^G_{k} \leq Q_{k\,\max},\quad \forall k \in \mathcal{G},\\
	&\alpha_k \geq 0, \quad \forall k \in \mathcal{G}.
	\end{align}
\end{subequations}
Finally, the voltages should satisfy the following \textit{Voltage Constraints}
\begin{subequations}\label{eq: voltage constraints}
\begin{align}
\label{eq: voltage bound 1}
&V_{k\,\min}\leq |V_k|\leq V_{k\,\max}, \quad \forall k \in \mathcal{N}, \\
\label{eq: voltage bound 2}
& |V_l-V_m|\le \Delta V_{lm}^{\max},
 \quad \forall (l,m) \in \mathcal{L}.
\end{align}
\end{subequations}
In \cite{Lavaei:2012}, the constraints in \eqref{eq: voltage bound 2} have been proven to be practically equivalent to the more classical bounds
\[
 \left | V_l(V_l-V_m)^*y_{lm}^*\right | \le S_{lm\, \max},\quad\forall(l,m)\in\mathcal{L},
\]
where $S_{lm\, \max}$ is the maximum apparent power flow which can path through the line $(l,m)\in\mathcal{L}$.

\subsection{Control and State Variables}
\label{sec:control-state}
In the power-systems literature, the  variables appearing in formulation above are usually divided into  two classes.
Indeed, already in \cite{Carpentier:1979aa}, the distinction between  control and state variables is explicitly made. \emph{Control variables} are those used by the network operator to set the operating condition of the network. \emph{State variables} are \emph{dependent} variables that represent the state of a power network, see also the recent survey \cite{Capitanescu:2016}


In particular, control and state variables are defined differently depending on the type of bus. In a generator bus $k\in\mathcal{G}$, usually referred to as \textit{PV bus} (see e.g.\ \cite[Remark 1]{Low:2014a}) mean active power $P^G_{k}$ of the generator, deployment factor $\alpha_k$ and magnitude $|V_k|$ of the complex bus voltage represent the control variables, 
while phase angle $\theta_k$ of bus voltage and generator reactive power $Q^G_{k}$ are the state variables. In a load bus, or \textit{PQ bus}, the active and reactive power of the load $P^L_{k},Q^L_{k}$ and the active and reactive power of the renewable generators $P^R_{k},Q^R_{k}$ are given (their values are known to the network operator)
while magnitude and phase angle of bus voltage $|V_k|,\theta_k$ are state variables. A node to which both a generator and loads are connected is to be considered as a generator bus. 

{\begin{remark}[Slack bus]
In power flow studies, usually bus $0$ is considered as \emph{reference} or \emph{slack} bus. The role of the slack bus is to balance the active and reactive power in the power grid. The slack bus must include a generator. In most cases, the voltage magnitude and phase are fixed at the slack bus, whereas active and reactive generator power are variables.  Therefore, in a slack bus there is no control variable, while active and reactive  generator power $P^G_k$ and $Q^G_k$ are state variables. For simplicity of notation, we do not include slack bus in the formulation,. However, its introduction is seamless, and slack bus is considered in the numerical simulation. 
\end{remark}
}

To emphasize the inherent difference between control and state variables, we introduce the notation
\begin{align}\label{eq: control variables}
\uu &\doteq \{P^G_1\,  \cdots \,P^G_{n_g},
|V_1| \cdots |V_{n_g}|, \alpha_1,\ldots,\alpha_n\}\\
\label{eq: state variables}
\xx_{\dd} & \doteq \{Q^G_1\,  \cdots \,Q^G_{n_g},
|V_{n_g+1}|,\ldots,|V_n|,\theta_1,\ldots,\theta_n\}
\end{align}
Note that in the notation above, we emphasize the fact that the state variables 
depend on the realization of the uncertainty.

We are now ready to provide a formal statement of the  \emph{robust} version of the optimal power flow problem. To this end, denote by $g(\uu,\xx_{\dd},\dd) = 0$
the uncertain equality constraints collected in \eqref{eq:unc_bal}, and by $h(\uu,\xx_{\dd},\dd) \le 0$ the uncertain inequalities collected in \eqref{eq:unc power constraints} and \eqref{eq: voltage constraints}.

\vspace{.1in}
\noindent{\bf Robust AC-OPF}
\begin{eqnarray} \label{eq: Robust AC-OPF} 
\underset{\uu}{\min} \
\ f(\uu)
\   \text{s.t.:}&
\text{for all } \dd \in\DD ,
\text{ there exist } \xx_{\dd} \text{ such that }\ \ \ \\ \nonumber
&\qquad g(\uu,\xx_{\dd},\dd) = 0  \text{ and } h(\uu,\xx_{\dd},\dd) \le 0.
\end{eqnarray}

}

In the above formulation of the robust OPF problem, seemingly first introduced in \cite{Zhang2011,Wada:2014}, the objective is to optimize the values of the nominal generated power, the bus voltage magnitude at generator node and the deployment vector so that i) the network operates safely for all values of the uncertainty and ii) the generation cost of the network is minimized.
In fact, if a solution to the problem above exists, we guarantee that for any admissible uncertainty, there exists a network state $\xx_{\dd}$ satisfying the operational constraints.

 Note that problem \eqref{eq: Robust AC-OPF} is computationally very hard.
 Indeed, even when no uncertainty is present, the AC-OPF problem is known to be non-convex, due to the presence of the nonlinear (quadratic) equality constraints \eqref{eq:unc_act_bal} and \eqref{eq:unc_reac_bal}, and the nonconvex line constraints \eqref{eq: voltage constraints}.
On top of that, these nonconvex constraints should be guaranteed over all the uncertainty set $\DD $, which is an infinite set (actually, uncountable). Hence, the optimization problem  \eqref{eq: Robust AC-OPF} is a nonlinear/nonconvex semi-infinite optimization problem.
 
We remark that taking a robust approach, i.e.\ enforcing the constraints in~\eqref{eq: Robust AC-OPF} for ``all'' possible values of the uncertain parameters, is in many cases excessive, and leads to  conservative results, with consequent degradation of the cost function (in our case, higher generation cost). 
Hence, in the remainder of this paper, we follow a chance-constrained approach, in which a probabilistic description of the uncertainty is assumed, and a solution is sought which is valid for the entire set of uncertainty except for a (small) subset having probability measure smaller than  a desired (small) risk level~$\varepsilon$. This approach is suitable for problems where ``occasional'' violation of constraints can be tolerated. One can argue that this is the case in power networks, since violation of line flow constraints does not necessarily lead to immediate line tripping. Rather, the line gradually heats up until a  critical condition is reached and only then the line is disconnected. 
Therefore, if line overload happens with low probability, this will not lead to line tripping nor it will damage the network. 
Also, in the robust case, one clearly needs to assume that the uncertainty set $\DD$ is bounded, and to know this bound. In a probabilistic setup, however, one can also consider cases where $\DD$ is unbounded, that is where the distribution of $\dd$ has infinite support, as e.g.\ the Gaussian case.

Formally, in the sequel, we assume that the uncertainty vector $\dd$ is random with possibly unbounded support $\DD$. Then, given a (small) risk level~$\varepsilon\in(0,1)$, the chance constrained version of the optimal power flow problem is stated as follows

\vspace{.1in}
\noindent
\textbf{Chance-constrained AC-OPF}
\begin{eqnarray}  \label{eq: OPF CCP }
\underset{\uu}{\min}
\ f(\uu)
\ \text{s.t.:}& 
\Pr\Big\{\dd \in\DD\text{ for which }
\nexists\, \xx_{\dd} \text{ such that } \ \ \ \ \ \  \\ \nonumber
& g(\uu,\xx_{\dd},\dd ) = 0  \text{ and } h(\uu,\xx_{\dd},\dd  ) \le 0 \Big\}\leq \varepsilon. \nonumber
\end{eqnarray}	

We remark that, from a numerical point of view, the above chance-constrained optimization problem is not tractable. On one hand,  the original problem \eqref{eq: Robust AC-OPF}, which does not have probabilistic constraints, is, as mentioned before, non-convex. On the other hand,  the presence of probabilistic constraints in \eqref{eq: OPF CCP } requires the solution of hard multi-dimensional integration problems.
\FD{To circumvent these numerical difficulties, in the next section we introduce an integrated solution approach that allows to i) relax the non-convex problem \eqref{eq: Robust AC-OPF} to a convex one and ii) derive randomized algorithms for solving the chance constraint optimization problem.}

\section{Efficient Numerical Relaxation}
\label{sec:relaxation}


\subsection{Approximation of Non-Convex Terms}
In this subsection, we extend the relaxation technique discussed in \cite{Bai:08,Lavaei:2012} to the robust case circumventing the non-convexity associated with the optimal power flow problem \eqref{eq: Robust AC-OPF}. We stress that the only source of non-convexity of \eqref{eq: Robust AC-OPF}, is due to non-linear terms $V_kV_l$'s appearing in \emph{Balance Equations} \eqref{eq:unc_bal}, and \emph{Voltage Constraints} \eqref{eq: voltage constraints}. However, the quadratic constraints can be reformulated as linear ones by introducing a new variable $\mathbf{W}=\mathbf{V}\mathbf{V}^*$ where $\mathbf{V}$ is the vector of complex bus  voltages $\mathbf{V} \doteq[V_1,\ldots,V_n]^T$. In order to replace $\mathbf{V}\mathbf{V}^*$ with the new variable $\mathbf{W}$, two additional constraints need to be included: i) the matrix $\mathbf{W}$ needs to be positive semi definite, i.e. the following \textit{Positivity Constraint} should hold
$\mathbf{W}\succeq 0,$
and ii) its rank should be one, i.e.\ the following \textit{Rank Constraint} should hold
$\text{rank}\{\mathbf{W}\}=1.$ 
An important observation is that, by introducing  matrix $\WW$, the only source of nonconvexity is captured by the rank constraint. Indeed, as shown first in \cite{Lavaei:2012} and subsequently in \cite{madani2015convex},  in most cases this constraint can be dropped without affecting the OPF solution. 

To formally define the \textit{convexified version} of the Robust AC-OPF, we note that bus voltage $\mathbf{V}$ appears in \emph{Voltage Constraints} \eqref{eq: voltage constraints} and \emph{Balance Equations} \eqref{eq:unc_bal}. Therefore, these constraints are  redefined in terms of the new variable $\mathbf{W}$ as
\begin{subequations}\label{eq: voltage constraints W}
	\begin{align}
	\label{eq: voltage bound W}
	& (V_{k\,\min})^2\leq W_{kk}\leq (V_{k\,\max})^2, \,\, \forall k \in \mathcal{N}\\
\label{eq: voltage bound 2W}
& W_{ll}+W_{mm}-W_{lm}-W_{ml}\le (\Delta V_{lm}^{\max})^2,\,\,\forall(l,m)\in\mathcal{L}
	\end{align}
\end{subequations}
\small \begin{subequations}\label{eq:unc_bal W}
	\begin{align}
	\nonumber
	&P^G_{k}+\alpha_k \mathbf{s}^T \real\{\dd\} + P^R_{k}(\dd) -  P^L_{k}(\dd)  =\\ \label{eq:unc_act_bal W} 
    & \qquad\sum_{l\in\mathcal{N}_k} \real\left\{(W_{kk}-W_{kl})^*y_{kl}^*\right\} \\
	\label{eq:unc_reac_bal W}
	&Q^G_{k}  +  Q^R_{k}(\dd) - Q^L_{k}(\dd) =
	 \sum_{l\in\mathcal{N}_k} \imag\left\{(W_{kk}-W_{kl})^*y_{kl}^*\right\},~\forall k \in \mathcal{N}
	\end{align}
\end{subequations}
\normalsize
respectively. 
Recalling that the first $n_g$ buses are generator buses and the remaining ones are load buses, we observe that the matrix $\mathbf{W}$ can be written as
\small 
\begin{equation}
\mathbf{W}=\mathbf{V}\mathbf{V}^* = \left[\begin{array}{cccccc}
|V_{1}|^2  & V_{1}V_{2}^* & \ldots  &  V_{1}V_{n_g}^* &  \ldots & V_{1}V_{n}^* \\
  &  |V_{2}|^2  & V_{2}V_{3}^*  &  \ldots   &  \ldots   & V_2V_n^*\\
&  &  \ddots &     &  &  \\

&  &   & |V_{n_g}|^2    &  &  \\
&  &   &   &   \ddots &  \\
&  &   &   &  &  |V_{n}|^2 \\
\end{array}\right].
\end{equation}
\normalsize
It is  immediately noticed that some elements of $\WW$ involve the control variables $|V_{1}|,\ldots,|V_{n_g}|$, 
 while  others are dependent variables corresponding to the voltage magnitude $|V_{n_g+1}|,\ldots,|V_n|$ at non-generator nodes, and the voltage phases $\theta_1,\ldots,\theta_n$. 
In order to distinguish between control and state variables appearing in it, we ``decompose'' $\WW$ into the sum of two submatrices $\mathbf{W}^\mathbf{u}$ and  $\mathbf{W}^\mathbf{x}$ as follows
\begin{equation} \label{eq: W_u}
\mathbf{W}^\mathbf{u} \doteq \text{diag}(|V_1|^2, \ldots,|V_{n_g}|^2, 0,\ldots),  
\mathbf{W}^\mathbf{x} \doteq \mathbf{W}-\mathbf{W}^\mathbf{u}.
\end{equation}
In this decomposition, we have a matrix $\WW^\mathbf{u}$ that includes the diagonal elements of $\WW$ corresponding to the generator nodes only, while the remaining elements of $\WW$ are collected in $\WW^\mathbf{x}$. With this in mind, the control and state variables are redefined as 
$\uu \doteq \{\mathbf{P}^{G}, \aa,\WW^\uu\}$
and 
$\xx_{\dd} \doteq \{\mathbf{Q}^G,\WW^\xx \}$
respectively. 
With this notation settled, we are in the position to formally introduce the convexified version of the robust AC optimal power flow problem as follows\\
\vskip -0.1in
\noindent{\bf Convexified Robust AC-OPF (CR-AC-OPF)}
{\small
\begin{align}  \label{eq: CR-AC-OPF}
&\underset{\mathbf{P}^{G}, \aa,\WW_\uu}{\text{minimize }} 
\sum_{k\in\mathcal{G}} f_k(P^G_{k})\\ \nonumber
&\text{subject to: } 
\text{for all } \dd \in\DD ,
\text{ there exist } \mathbf{Q}^G=\mathbf{Q}^G(\dd),
\WW^\xx=\WW^\xx(\dd) \\\nonumber
&\text{such that } \\ \nonumber
&\quad \nonumber
\WW=\WW^\uu+\WW^\xx, 
\quad \mathbf{1}^T\aa=1,  \quad \WW\succeq 0,\quad \alpha_k\geq0,  \quad \forall k\in\mathcal{G} \\
&\quad  \nonumber
P^G_{k}+\alpha_k \mathbf{s}^T \real\{\dd\} + P^R_{k}(\dd) -  P^L_{k}(\dd) =\\ \nonumber &\qquad\qquad\qquad\sum_{l\in\mathcal{N}_k} \real\left\{(W_{kk}-W_{kl})^*y_{kl}^*\right\},  \; \forall k \in \mathcal{N}\\
&\quad  \nonumber
Q^G_{k} +  Q^R_{k}(\dd) - Q^L_{k}(\dd) =
\sum_{l\in\mathcal{N}_k} \imag\left\{(W_{kk}-W_{kl})^*y_{kl}^*\right\}, \forall k \in \mathcal{N}\\
&\quad  \nonumber
P_{k\,\min}\leq P^G_{k}+\alpha_k \mathbf{s}^T \real\{\dd\} \leq P_{k\,\max}, \quad \forall k \in \mathcal{G} \\
&\quad  \nonumber
Q_{k\,\min}\leq Q^G_{k}\leq Q_{k\,\max}, \quad \forall k \in \mathcal{G}\\
&\quad  \nonumber
\left(V_{k\,\min}\right)^2\leq W_{kk}\leq \left(V_{k\,\max}\right)^2,  \quad \forall k \in \mathcal{N}\\
&\quad  \nonumber
W_{ll}+W_{mm}-W_{lm}-W_{ml}\le (\Delta V_{lm}^{\max})^2,\ \quad\forall(l,m)\in\mathcal{L}.
\end{align}} \vskip -.2in

\begin{remark}[About \textbf{CR-AC-OPF} formulation]
\label{crac-remark}
To the best of the authors knowledge, the formulation of the  \textbf{CR-AC-OPF} is original, and it represents the first main contribution of the present paper. First, the approach clearly differs from the formulations based on DC power flow. Second, it improves upon the formulation based on convex AC-formulation  in \cite{ETH-AC}.

Indeed, in that work, the authors cope with the need of guaranteeing the existence of a different value of $\WW$ for different values of the uncertainty $\dd$ by imposing a specific dependence on $\WW$ from the uncertainty. In our notation, \cite{ETH-AC} introduces the following finite (linear) parameterization
\begin{equation}
\label{eq:W-param}
\WW(\dd)=A+\sum_k B_k\dd_k,
\end{equation}
where $A, B_1,\ldots,B_n$ become design variables in the optimization problem. 
{

We note that the choice of introducing some type of parameterization of the voltages with respect to uncertainty is unavoidable if a certificates approach is not adopted. this is due to the equality constraints in the Balance Equations. 
In DC-based results, these become linear equalities, and hence they can be explicitly solved, as in as \cite{Bienstock2014,Lubin2016}.
In nonlinear approaches, this parameterization is always found, and goes under different names;
for instance it is referred to as \textit{affine control}
\cite{Morari2015}, 
\textit{affinely adjustable robust counterpart} in \cite{Jabr2013}, \textit{affine feedback policies} in \cite{Guo2018}, and first-order Taylor approximation in \cite{Andersson2018}.

With respect to these formulation, we remark the following:
i) the \textbf{CR-AC-OPF} formulation is surely less conservative, since it does not impose a \textit{specific dependence} on $\WW$. The gain in terms of performance with respect to \cite{ETH-AC} is reported in the numerical experiments in  Section \ref{sec:examples}, 
ii)
the parameterization in \eqref{eq:W-param} requires modification of   the voltage magnitude at the generators during operation. Indeed, as previously noted, $\WW$ involves not only the state variables,  but also the control variables corresponding to generator voltages $|V_{1}|,\ldots,|V_{n_g}|$ which are not supposed to change during the operation. In the \textbf{CR-AC-OPF}, contrary, all control variables---including generators voltage magnitude---are designed during the optimization phase and do not need to be changed when demand/renewable power fluctuations occur.}
\end{remark}

Moreover, although the \textbf{CR-AC-OPF} formulation represents a relaxation of the robust optimal power flow problem, we can show that there are special cases in which this relaxation is indeed exact. In particular, we introduce next a class of networks where \textbf{CR-AC-OPF} does provide the optimal solution to the Robust AC-OPF problem. To this end, we first briefly recall some concepts from graph theory that are central for the results to follow. For a more detailed discussion on the following definitions and their interpretation in the context of power networks the reader is referred to \cite{madani2015convex}.

\begin{definition}
	A network/graph is called \emph{weakly cyclic} if every edge belongs to at most one cycle.
\end{definition}
\begin{definition}
	A network is said to be \emph{lossless} if
    \mbox{$
    \real\{y_{lm}\} = 0$,} $ \forall (l,m) \in \mathcal{L}.
    $
\end{definition}

With these definitions at hand, we can now describe a  class of networks for which the relaxed formulation above provides an exact solution. This theorem represents  a natural extension of the results presented in \cite{madani2015convex} to the robust case and the proof is given in appendix.
\begin{theorem}
	\label{theo-OPF-robust}
	Consider a lossless weakly-cyclic network with cycles of size 3, and assume $Q_k^{\min}=-\infty$ for every $k\in\mathcal{G}$.  Then, the convex relaxation {\bf CR-AC-OPF} is exact.
\end{theorem}
\noindent
\textbf{Proof:} See Appendix \ref{sec: appendix}.

\begin{remark}
In the theorem above, by exact we mean that i) the robust convex relaxation has an optimal value equal to the one of the  Robust AC-OPF and ii) for any admissible value of the uncertainty $\dd$, there exists a rank one matrix $\WW_\delta$ that satisfies the constraints of the convex relaxation and, hence, there exists a $\VV_\delta$ that satisfies the constraints of the Robust AC-OPF.
\FD{Hence,  if the convex relaxation \textbf{CR-AC-OPF} is feasible, then  Theorem~1 guarantees that the ensuing optimal values of $\mathbf{P}^{G}$, $\aa$,$\mathbf{W}_\mathbf{u}$ are such that for every value of the uncertainty there exists a rank one $\mathbf{W}$ solving the problem. 
Hence, for those cases when a solution to \textbf{CR-AC-OPF} is found, the Theorem provides a very strong guarantee.
On the other hand, it should be remarked that the theoretical guarantees of Theorem 1 clearly do not automatically translate for those  cases when the robust solution of \textbf{CR-AC-OPF} is {prohibitive from a computational viewpoint}.
In these cases, the existence of a rank one solution for {some uncertainty instances} may not imply the existence of a rank one solution for an unseen {instance of the uncertainty}. 
}\end{remark}

We should remark that the class considered in Theorem~\ref{theo-OPF-robust} is not fully realistic. However, 
we feel that this result is important for several reasons. First, it shows that  the exactness of the relaxation for a particular class of networks---proven for the nominal case in \cite{madani2015convex}---carries over to the robust formulation we introduced in \eqref{eq: CR-AC-OPF}. This shows that {\bf CR-AC-OPF} represents indeed the right way to formulate the robust counterpart of the AC-OPF problem. Second, 
several results of the same spirit have been obtained in the literature for deterministic (no uncertainty) networks, see for instance \cite{Low:2014b}: we believe that these results can be extended to the convexified robust formulation introduced here. This will the subject of further research.

\subsection{Penalized and robust costs}\label{sec: worst case cost}
Note that, besides the  configurations where the relaxation has been proven to be exact, in the general case it is important to obtain solutions $\WW$ with low rank. To this end, we adopt a penalization technique similar to that introduced in  \cite{Madani_promises,madani2015convex}. More precisely, it is observed there that maximizing the weighted sum of off-diagonal entries of $\WW$ often results in a low-rank solution. 
To this end, one can augment the objective function with a weighted sum of generators' reactive power---for lossless network---and apparent power loss over the series impedance of some of the lines of the network (the so-called problematic lines $\mathcal{L}_{\text{prob}}$)---for lossy network---to increase the weighted sum of the real parts of off-diagonal elements of $\WW$ and hence obtain a low rank solution. 
Formally, let $L_{lm} \doteq|S_{lm}+S_{ml}|$  denote the apparent power loss over the line $(l,m)$ with $S_{lm}\doteq|(W_{ll}-W_{lm})^*y_{lm}^*|$. Given nonnegative  factors $\gamma_b$ and $\gamma_\ell$, in  \cite{Madani_promises,madani2015convex} the following  penalized cost was considered
\[
f_\mathrm{pen}(\mathbf{P}^G,\mathbf{Q}^G,\WW)\doteq
\sum_{k\in\mathcal{G}} f_k(P^G_{k}\!)+ \gamma_b \!\!\sum_{k\in\mathcal{G}}\!\! Q^G_{k}+ \!\gamma_\ell  \!\!\!\!\sum_{(l,m)\in\mathcal{L}_{\text{prob}}}\!\!   \!\!\!\!\!\!L_{lm}.
\]
Since, in the proposed approach,  $\mathbf{Q}^G$ and $\WW$ depend on the value of the uncertainty $\dd$, the formulation of the robust optimal power flow problem needs to reflect this fact. More precisely, the right cost function to consider is
 \[
\max_{\dd \in \DD} f_\mathrm{pen}(\mathbf{P}^G,\mathbf{Q}^G,\WW)\doteq
\]
\[
 \max_{\dd \in \DD} \left(
\sum_{k\in\mathcal{G}} f_k(P^G_{k}\!)+ \gamma_b \!\!\sum_{k\in\mathcal{G}}\!\! Q^G_{k}+ \!\gamma_\ell  \!\!\!\!\sum_{(l,m)\in\mathcal{L}_{\text{prob}}}\!\!   \!\!\!\!\!\!L_{lm}\right).
\]
The constraints remain the same as in~\eqref{eq: CR-AC-OPF}. This modification of the cost provides a way of ``encouraging'' low rank solutions and it has been shown to work well in practice. 

\section{A randomized approach to  \textbf{CR-AC-OPF}}
\label{sec:scenario}
In this section, we first briefly summarize  the main features of the scenario with certificates problem, and then show how this approach can be used to tackle in an efficient way the \textbf{CR-AC-OPF} problem introduced in the previous section.

\subsection{Scenario with Certificates}

In \cite{SWC-TAC} a \FD{class of robust optimization problems ``with certificates" was formally introduced. This definition refers to robust optimization problems in which a clear distinction can be made between so-called \textit{design} variables $\theta$ and \textit{certificates} $\xi$. 
While the design variables  $\theta$ play the role of classical optimization variable, the certificates are variables for which we are not interested in, as long as we are guaranteed that for any possible value of the uncertainty there exists a value of $\xi$ guaranteeing feasibility of the problem.}

\FD{Formally, we assume that the design variables $\theta$ belong to a given convex set $\Theta\subseteq\Real{n_\theta}$,
and the certificates $\xi$ belong to a compact (possibly nonconvex) set $X\subseteq\Real{n_\xi}$. Then, we
consider a function $f(\theta,\xi,\dd)$ which is jointly convex in $\theta$ and $\xi$, for any given $\dd\in\DD$. Then, we can define the so-called  \textit{robust optimization problem with certificates}  introduced as follows} 
\begin{eqnarray}
\min_{\theta} && c^T \theta \label{eq:certificates_opt}\\
                  \text{subject to} && \forall \dd\in\DD  \ \exists \xi=\xi(\dd)
\text{ satisfying } f(\theta,\xi,\dd)\le 0. \nonumber
\end{eqnarray}
Then, it is shown that problem \eqref{eq:certificates_opt} can be approximated by introducing an appropriate randomized counterpart, based  on the extraction of $N$ random samples of the uncertainty.
Formally, we extract an \textit{uncertainty multisample}
\[
\{\dd^{(1)},\ldots,\dd^{(N)}\}
\]
and construct the following \textit{scenario problem with certificates}
\begin{eqnarray}
         \tSOC=\arg\min_{x,\xi_{1},\ldots,\xi_{N}}&& c^T \theta  \label{eq:scenario_cert}\\
                      \text{subject to}&&  f(\theta,\xi_{i},\dd^{(i)}) \leq 0, \  i=1,\ldots,N.\nonumber
\end{eqnarray}
Note that 
contrary to the classical scenario problem, 
in SwC \textit{a new certificate variable $\xi_{i}$ is created for every sample $\dd^{(i)}$}. In this way, one implicitly constructs an "uncertainty dependent" certificate, without assuming any \textit{a-priori} explicit functional dependence on $\dd$.

To present the properties of the solution $\tSOC$, we first introduce the violation probability of a given design $\theta$ as follows
\[
\V(\theta)=
\Pr\Bigl\{\exists \dd\in\DD  | \nexists \xi \text{ satisfying } f(\theta,\xi,\dd)\le 0 \Bigr\}.
\]
Then, the main result regarding the scenario optimization with certificates is recalled next. The result was derived in \cite{SWC-TAC} as an extension to the classical scenario approach developed in~\cite{calafiore2006scenario}.
\vskip .2in
\begin{theorem}
\label{them:SwC}
Assume that, for any multisample extraction, problem \eqref{eq:scenario_cert} is feasible and attains a unique optimal solution.
Then, given an accuracy level $\varepsilon\in(0,1)$ and a confidence level $\beta\in(0,1)$, if the number of samples is chosen as
\begin{equation}
\label{eq: N_SwC}
N\ge N\ped{SwC} = \frac{\mathrm{e}}{\varepsilon(\mathrm{e}-1)}\left(\ln\frac{1}{\beta}+n_{\theta}-1\right)
\end{equation}
where $n_\theta$ is the dimension of $\theta$, and e is the Euler number. Then, with probability at least $1-\beta$ , the solution $\tSOC$ of problem \eqref{eq:scenario_cert} satisfies
$
\V(\tSOC)\le\varepsilon.
$
\end{theorem}

Note that, in the above theorem, we guaranteed with high confidence $(1-\beta)$ that the solution returned by the SwC has a risk of violation of constraints less than the predefined (small) risk level~$\varepsilon$.
\subsection{SwC Solution to \textbf{CR-AC-OPF}}
Clearly, problem \textbf{CR-AC-OPF} represents a 
robust optimization problem with certificates, in which the design variables are those "controllable" by the network manager, i.e.
$\theta\equiv \uu$
while the certificates represent the quantities that can be "adjusted" to guarantee constraint satisfaction, i.e.
$\xi \equiv \xx_{\dd}$.
Based on this consideration, in this paper we propose the following solution strategy to the \textbf{CR-AC-OPF} problem: \\

\noindent{\bf SwC-AC-OPF Optimization Procedure}
\begin{itemize}
\item[i)] Given probabilistic levels $\varepsilon$, and $\beta$ compute $N\ped{SwC}$ according to \eqref{eq: N_SwC}.
\item[ii)] Generate $N\ge N\ped{SwC}$ sampled scenarios $\dd^{(1)},\ldots,\dd^{(N)}$, where the uncertainty is drawn according to its known probability density.
\item[iii)]
Solve the following convex optimization problem, which returns the control variables $\mathbf{P}^{G}, \mathbf{W^u},  \aa$.\\

\noindent{\bf SwC-AC-OPF} {\small
\begin{align}  \label{eq: SwC Original}
&\underset{\tiny 
\mathbf{P}^{G}, \mathbf{W^u},  \aa,
\mathbf{Q}^{G,[1]},\ldots,\mathbf{Q}^{G,[N]},
\WW^{\xx,[1]},\ldots, \WW^{\xx,[N]}
}{\text{minimize}} 
\ {\boldsymbol{\gamma}}\\
 &\text{subject to:} \text{ for } i=1,\ldots,N\nonumber \\
&\nonumber \WW^{[i]}=\WW^\uu+\WW^{\xx,[i]}, \quad \WW^{[i]}\succeq 0, \quad \alpha_k\geq0.  \quad \forall k\in\mathcal{G}\\
&\nonumber
L_{lm}^{[i]} =
|(W^{[i]}_{ll}-W^{[i]}_{lm})^*y_{lm}^*|
+
|(W^{[i]}_{mm}-W^{[i]}_{ml})^*y_{lm}^*|\\
&\sum_{k\in\mathcal{G}} f_k(P^G_{k}\!)+ \gamma_b \!\!\sum_{k\in\mathcal{G}}\!\! Q^{G,[i]}_{k}+ \!\gamma_\ell  \!\!\!\!\sum_{(l,m)\in\mathcal{L}^{\text{prob}}}\!\!   \!\!\!\!\!\!L_{lm}^{[i]} \le \gamma
\nonumber\\
&\nonumber
P^G_{k}+\alpha_k \mathbf{s}^T \real\{\dd^{(i)}\} + P^R_{k}(\dd^{(i)}) -  P^L_{k}(\dd^{(i)}) = \\ \nonumber
& \qquad\qquad\qquad\sum_{l\in\mathcal{N}_k} \real\left\{(W_{kk}^{[i]}-W_{kl}^{[i]})^*y_{kl}^*\right\},  \quad \forall k \in \mathcal{N}\\
&\nonumber
Q^{G,[i]}_{k} +  Q^R_{k}(\dd^{(i)}) - Q^L_{k}(\dd^{(i)}) = \\ \nonumber
& \qquad\qquad\qquad\sum_{l\in\mathcal{N}_k} \imag\left\{(W_{kk}^{[i]}-W_{kl}^{[i]})^*y_{kl}^*\right\}, \quad \forall k \in \mathcal{N}\\
&\nonumber
P_{k\,\min}\leq P^G_{k}+\alpha_k \mathbf{s}^T \real\{\dd^{(i)}\} \leq P_{k\,\max}, \quad \forall k \in \mathcal{G} \\
&\nonumber 
Q_{k\,\min}\leq Q^{G,[i]}_{k}\leq Q_{k\,\max}, \quad \forall k \in \mathcal{G}\\
&\nonumber 
\left(V_{k\,\min}\right)^2\leq W_{kk}^{[i]}\leq \left(V_{k\,\max}\right)^2,  \quad \forall k \in \mathcal{N}\\
&\nonumber
W_{ll}^{[i]}+W_{mm}^{[i]}-W_{lm}^{[i]}-W_{ml}^{[i]}\le (\Delta V_{lm}^{\max})^2,\ \quad\forall(l,m)\in\mathcal{L} 
\end{align}}

\item[iv)] During operation, measure  uncertainty in generations and loads $\dd$, and accommodate the $k$-th controllable generator as 
$\bar{P}^G_{k} =P^G_{k}+\alpha_k \mathbf{s}^T \real\{\dd\},~ 
|V_{k}|  =\sqrt{W_{kk}},~k \in \mathcal{G}$. 
\end{itemize}

We remark that Theorem~\ref{them:SwC} guarantees that the \textbf{SwC-AC-OPF} optimization procedure is such that the probabilistic constraints of the Chance-Constrained AC-OPF \eqref{eq: OPF CCP } are satisfied with high confidence~$(1-\beta)$.
In other words, one has an \textit{a-priori} guarantee that the risk of violating the constrains is bounded, and one can accurately bound this violation level by choosing the probabilistic parameters $\varepsilon,\beta$. The above described procedure represents the main result of this paper. In Section \ref{sec:examples}, we demonstrate how this procedure outperforms existing ones in terms of guarantees of lower line-violations.

\begin{remark}[Using available data on uncertainty]
\FD{It is frequently the case that several measurements of the uncertainty are available to the system designer.
Several approaches aiming at designing uncertainty sets based on this available data have been hence recently proposed in the literature, and employed in approaches such as Robust Sample Average Approximation and Distributionally Robust methods, see \cite{Guo2018} and references therein.
We remark that it is not directly necessary in a sample-based approach---as the one presented here---to design the uncertainty set. Indeed, the data itself can be used as ``scenarios" to be fed to the SwC problem, thus  avoiding this critical point. }
\end{remark}

\subsection{Handling Security Constraints}
An important observation is that the approach introduced in this paper can be readily extended to handle security constraint. In particular, the popular $N-1$ security requirements discussed, for instance, in \cite{ETH-AC} can be directly included in the robust optimization problem \eqref{eq: Robust AC-OPF}, and in its subsequent derivations.
We recall that, in the $N-1$ security constrained OPF framework, only the outages of \textit{a single component} are taken into account. A list of $N_\mathrm{out}$ possible outages, $\mathcal{I}^\mathrm{out}=\left\{0,1,\ldots,N_\mathrm{out}\right\}$, is formed (with $0$ corresponding to the case of no outages). Then, a 
large optimization problem is constructed, with $N_\mathrm{out}$ instances of the constraints, where the $i$-th instance corresponds to removing the $i$-th component form the equations, with $i\in\mathcal{I}^\mathrm{out}$. 
The exact same approach can be replicated here. 

However, we remark that the SwC formulation allows for a more precise and realistic handling of the possible outages. Indeed, we note that in a real network some buses may have a larger probability of incurring into an outage, for instance because they are located in a specific geographical position, or because they employ older technologies.
To reflect this scenario, we associate to each element $i\in\mathcal{I}^\mathrm{out}$ in the network a given \textit{probability of outage} 
\[
p_i^\mathrm{out}\in\left[0,1\right],
\]
which can be different for each component,
and is supposed to be known to the network manager.
Large  values of $p_i$ correspond to  large probability of an outage occurring in the $i$-th component.
Then, in our framework, outages can be treated in the same way as uncertainties.
Formally, we can associate to each component a random variable 
$\delta^\mathrm{out}_i,$ $ i\in\mathcal{I}^\mathrm{out}$,
with Bernoulli density with mean $1-p_i^\mathrm{out}$. Then, the scenario with certificates approach is applied considering the extended uncertainty 
$\left\{\dd,\dd^\mathrm{out}\right\}$.

In practice, when constructing the \textbf{SwC-AC-OPF} problem, for each sampling instance the $i$-th component  is removed with probability $p_i^\mathrm{out}$. It is important to note that this approach goes beyond the standard $N-1$ security constrained setup, since:  i) \textit{multiple} simultaneous outages are automatically taken into account, ii) it allows weighting differently the different component.

\section{Numerical Examples} \label{sec:examples}

In order to examine the effectiveness of the proposed method, we perform extensive simulations using New England 39-bus system case adopted from \cite{matpower2011}. The network has $39$ buses, $46$ lines and  $10$ conventional generation units, and it is modified to include $4$  wind generators connected to buses $5$, $6$, $14$ and $17$. The renewable energy generators and all loads connected to different buses are  considered to be uncertain. In total, there are $46$ uncertain parameters in the network. The goal is to design  active power and voltage amplitude of all controllable generators as well as the distribution vector $\boldsymbol{\alpha}$ such that the generation cost is minimized while all the constraint of the network i.e. line flow, bus voltage and generators output constraint are respected with high probability.

\subsection{Numerical Results}
\begin{figure}
		\centering
		\begin{subfigure}{.48\textwidth}
		\centering
		\includegraphics[trim={0cm 0cm .1cm 0cm},clip,width=0.9\columnwidth]{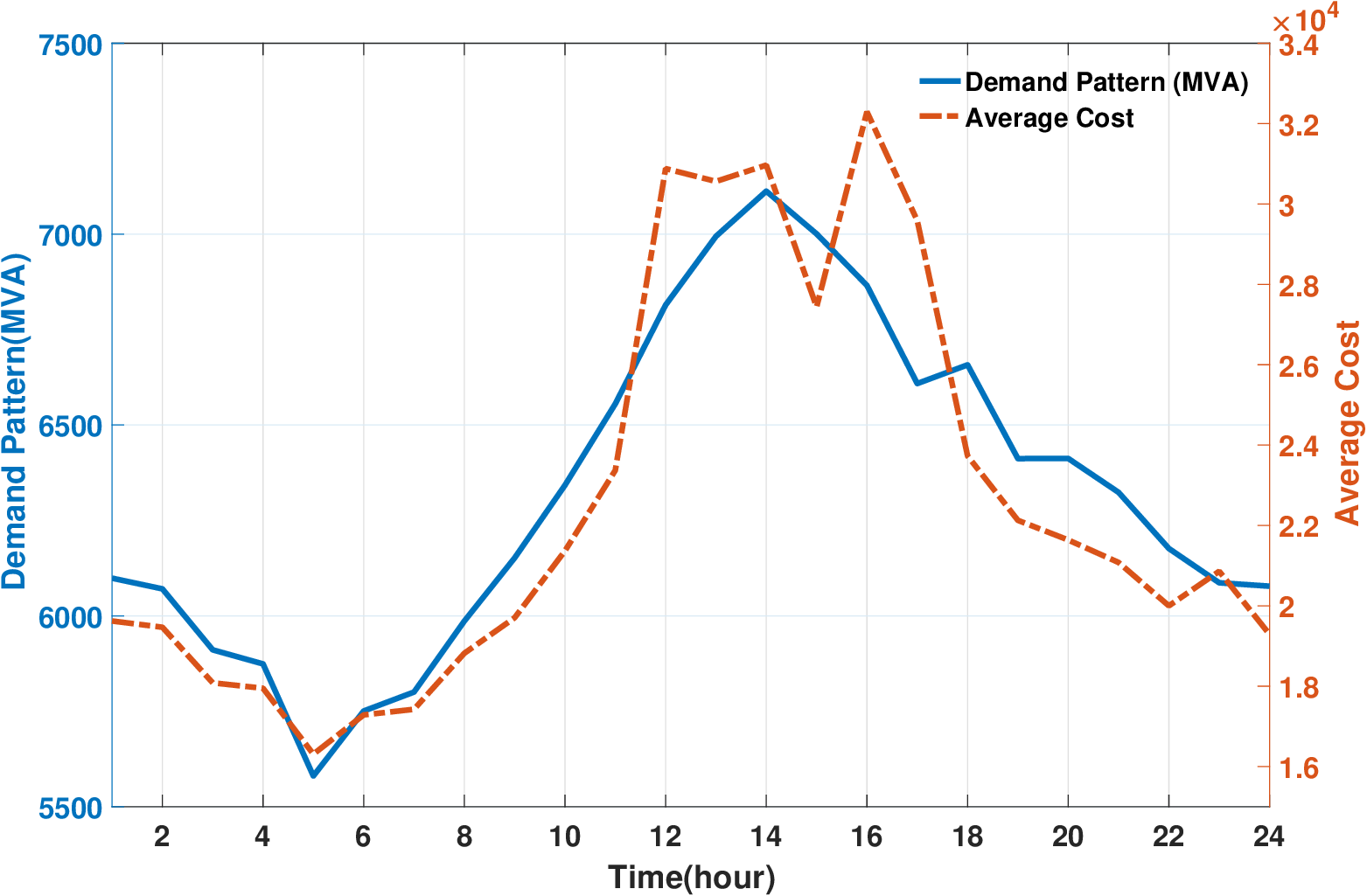}
\caption{}
\end{subfigure}
\begin{subfigure}{.5\textwidth}
		\centering
		\includegraphics[width=1.1\columnwidth]{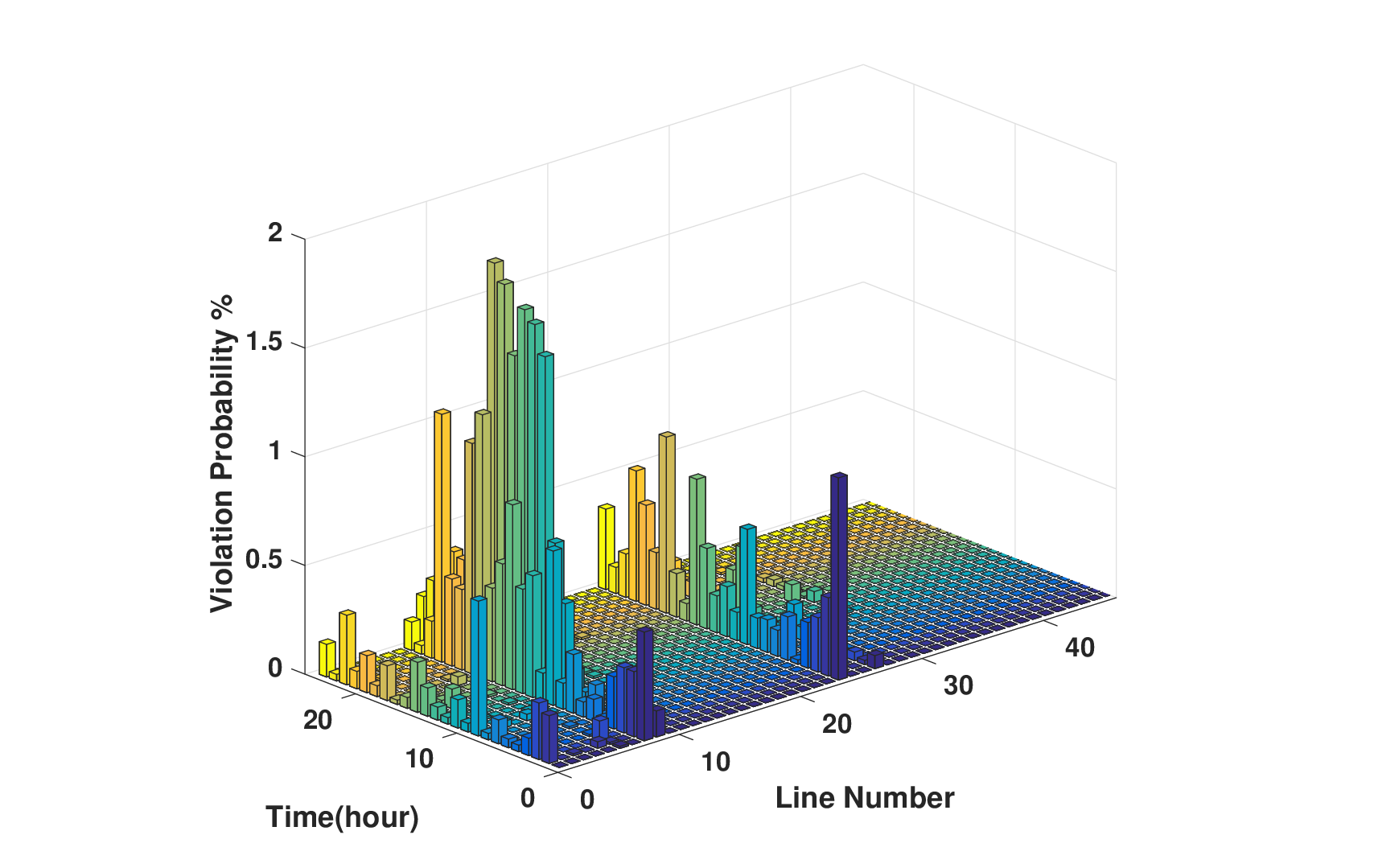}
\end{subfigure}
\caption{\footnotesize (a): $24$ hours demand pattern and average generation cost computed in the posteriori analysis based on Monte Carlo simulation. (b): $24$ hours empirical violation of line flow computed in the posteriori analysis based on Monte Carlo simulation. }
	\label{fig: 24hour results}
\end{figure}

We use the methodology presented in Section \ref{sec:scenario} to solve the robust optimal power flow problem for New England 39 bus system.
We consider a $24$ hour demand pattern shown in Fig. \ref{fig: 24hour results}-(a) and solve the optimal power flow problem for each hour by minimizing the nominal cost. The error probability distribution for wind generator and load is chosen based on Pearson system \cite{pearson1895contributions}---as suggested in \cite{hodge2013short,hodge_wind_2012}. Pearson system is represented by the  mean $\mu$ (first moment),  variance~$\sigma^2$ (second moment),  skewness $\gamma$ (third moment)  and 
kurtosis $\kappa$ (forth moment). In the simulation, we set $\sigma=0.2\times\text{(predicted value)},~\gamma=0,$ and, $\kappa=3.5$ leading to a  leptokurtic distribution with heavier tail than Gaussian. The selected probabilistic accuracy and confidence levels  $\varepsilon$ and $\beta$ are  $0.02$ and $1\times10^{-15}$ respectively resulting in $5,105$ number of scenario samples. We note that the number of design (control) variables  is $31$ in our formulation. 
The penetration level is chosen to be $30\%$.  
The penetration level indicates how much of the total demand is provided by the wind generators.
A $30\%$ penetration level means that, in total, the wind generators provide $30\%$ of the nominal load. We also assume that each renewable generator contributes equally to provide this power. 
The optimization problem is formulated and then solved using YALMIP \cite{loefberg_yalmip_2004} and Mosek \cite{andersen2000mosek} respectively, and returned the control variables $\mathbf{P}^{G}$, $\mathbf{W^u}$ and    $\aa$.

In order to examine  robustness of the obtained solutions, we run an  \textit{a-posteriori} analysis based on Monte Carlo simulation. To this end, we generated $10,000$ random samples from the uncertainty set---corresponding to uncertain active and reactive generated wind power and load---and for each sample, modified the network to include wind power generators, replaced the load vector by its uncertain counterpart, computed the power mismatch, distributed the mismatch to all conventional generators by using  (\ref{eq: affine control}) and solved the power flow problem whose constraint sets include the non-linear power balance equation, bus voltage constraints and generator power constraints. For each sample of the Monte-Carlo simulation, we solved the non-linear non-convex  feasibility problem to see if there exist feasible state variables that respect the constraints. We then  measured the number of infeasible samples as a measure of robustness of the designed control variables: the smaller the number of infeasible samples, the more robust are the control variables. This simulation is formally formulated in Algorithm \ref{alg: MonteCarlo}  where the nonlinear feasibility problems were solved using the function \texttt{fmincon} of Matlab\textregistered. 
It should be noted that the matrices $\mathbf{W}$ obtained were mostly not rank 1, but as it can be seen in the simulation, the control variable values obtained work well in practice.

Running Algorithm \ref{alg: MonteCarlo} with the control variables designed using the proposed SwC method results in only $14$ out of $10,000$ infeasible samples while the nominal control variables---the control variables designed without considering uncertainty---results in $5,374$ infeasible samples.   Hence, the probability of
joint feasibility of all constraints is simply approximated by [number of feasible instances]$/N$. In our experiment this probability turned out to be $9,986/10,0000=0.9986$, with a consequent probability of failure less than $.002$ (and hence much less than the considered $\epsilon=0.02$). 
On the contrary, the same analysis on the nominal design led to a joint violation  probability of $5374/10000=.5374$, hence more than $50\%$. The empirical violation is  shown in Table \ref{tab: comarison of empirical violation}.  We remark that both control variables are designed for the peak demand level which is at hour $14$. This  proves that the control variables designed using SwC method are indeed robust compared to design variables for which no uncertainty has been considered.

\renewcommand{\algorithmicrequire}{\textbf{Input:}}
\renewcommand{\algorithmicensure}{\textbf{Output:}}
\begin{algorithm}[H]
\begin{algorithmic}[1] \begin{footnotesize}
	\caption{Monte-Carlo Simulation}
	\label{alg: MonteCarlo}
	\STATE\algorithmicrequire{ Control variables $\{P^G_1\,  \cdots \,P^G_{n_g},
		|V_1| \cdots |V_{n_g}|, \alpha_1,\ldots,\alpha_n\}$}
	\STATE\algorithmicensure{\texttt{Infeas\_Counter}}\\
{\bf Initialization:}
\STATE Set \texttt{Infeas\_Counter = 0}\\
{\bf Evolution:}\\
Extract $10,000$ i.i.d samples $\dd^{(1)},\ldots,\dd^{(10,000)}$\\ 
\FOR{$i = 1$ \TO $10,000$}
\STATE{ Solve
 \begin{align*}
\text{find } \xx_{\dd}& \text{ s.t.}\\
&\nonumber
P^G_{k}+\alpha_k \mathbf{s}^T \real\{\dd^{(i)}\} + P^R_{k}(\dd^{(i)}) -  P^L_{k}(\dd^{(i)}) = \\ \nonumber
& \qquad\qquad\qquad\sum_{l\in\mathcal{N}_k} \real\left\{V_k(V_k-V_l)^*y_{kl}^*\right\},  \quad \forall k \in \mathcal{N}\\
&\nonumber
Q^{G}_{k} +  Q^R_{k}(\dd^{(i)}) - Q^L_{k}(\dd^{(i)}) = \\ \nonumber
& \qquad\qquad\qquad\sum_{l\in\mathcal{N}_k} \imag\left\{V_k(V_k-V_l)^*y_{kl}^*\right\}, \quad \forall k \in \mathcal{N}\\
&\nonumber
P_{k\,\min}\leq P^G_{k}+\alpha_k \mathbf{s}^T \real\{\dd^{(i)}\} \leq P_{k\,\max}, \quad \forall k \in \mathcal{G} \\
&\nonumber 
Q_{k\,\min}\leq Q^{G}_{k}\leq Q_{k\,\max}, \quad \forall k \in \mathcal{G}\\
&\nonumber 
V_{k\,\min}\leq V_{k}\leq V_{k\,\max},  \quad \forall k \in \mathcal{N}\\
&\nonumber
|V_l(V_l-V_m)^*y_{lm}^*|\leq S_{lm \max}, \quad \forall(l,m)\in\mathcal{L}
\end{align*}
\IF{$\xx_{\dd}==\emptyset$} 
\STATE{\texttt{Infeas\_Counter = Infeas\_Counter+1}}\\
\ENDIF
}
\ENDFOR \end{footnotesize}
 \end{algorithmic}
\end{algorithm}

\begin{table}[t]
\caption{Empirical violation of the Monte Carlo Algorithm \ref{alg: MonteCarlo} for the SwC solution and the nominal one---where no uncertainty is considered in the design phase.}
\begin{center}
\begin{tabular}{c||c}
\toprule
Approach & Empirical Violation\tabularnewline
\midrule
SwC & $0.0014$ \tabularnewline
\midrule
Nominal & $0.5374$ \tabularnewline
\bottomrule
\end{tabular}
\end{center}
\label{tab: comarison of empirical violation}
\end{table}

The rate of violation of network constraints can not be identified by the Monte Carlo simulation presented in Algorithm \ref{alg: MonteCarlo}. An infeasible optimization problem does not define if the infeasibility is due to violation of bus voltages, line ratings or generator active and reactive generated power. Furthermore, even if the constraint(s) causing infeasibility becomes clear, it remains unclear how much the constraint(s) exceed their limits. For this reason, we used \texttt{runpf} command of MATPOWER \cite{matpower2011} to derive line flows. We remark that Matpower power flow solver does not respect voltage and generator power constraints and only solves the non-linear balance equation. On the other hand---due to non-linearity of balance equation---the solution to power flow equation may not be unique \cite{korsak1972question,johnson1977extraneous,wang2003existence}. Hence, one may not obtain a solution that respects voltage and power constraints. For this reason we focus on line flows derived from MATPOWER.
The empirical violation of line flows is computed by counting the number of times each line exceeds its limit in the Monte Carlo simulation and dividing this value by $10,000$. Figure \ref{fig: 24hour results}(b) shows the $24$ hour empirical violation of line flow for all lines of the network. The mean value of generation cost associated with the $24$ hour demand is also shown in Fig. \ref{fig: 24hour results}(a). We note that this graph is obtained by averaging the generation cost over all $10,000$ samples. 

\begin{figure*}[t]
		\centering
	\begin{subfigure}{.48\textwidth}
		\centering
		\includegraphics[trim={.3cm 1.1cm 1.2cm 1cm},clip,width=0.9\columnwidth]{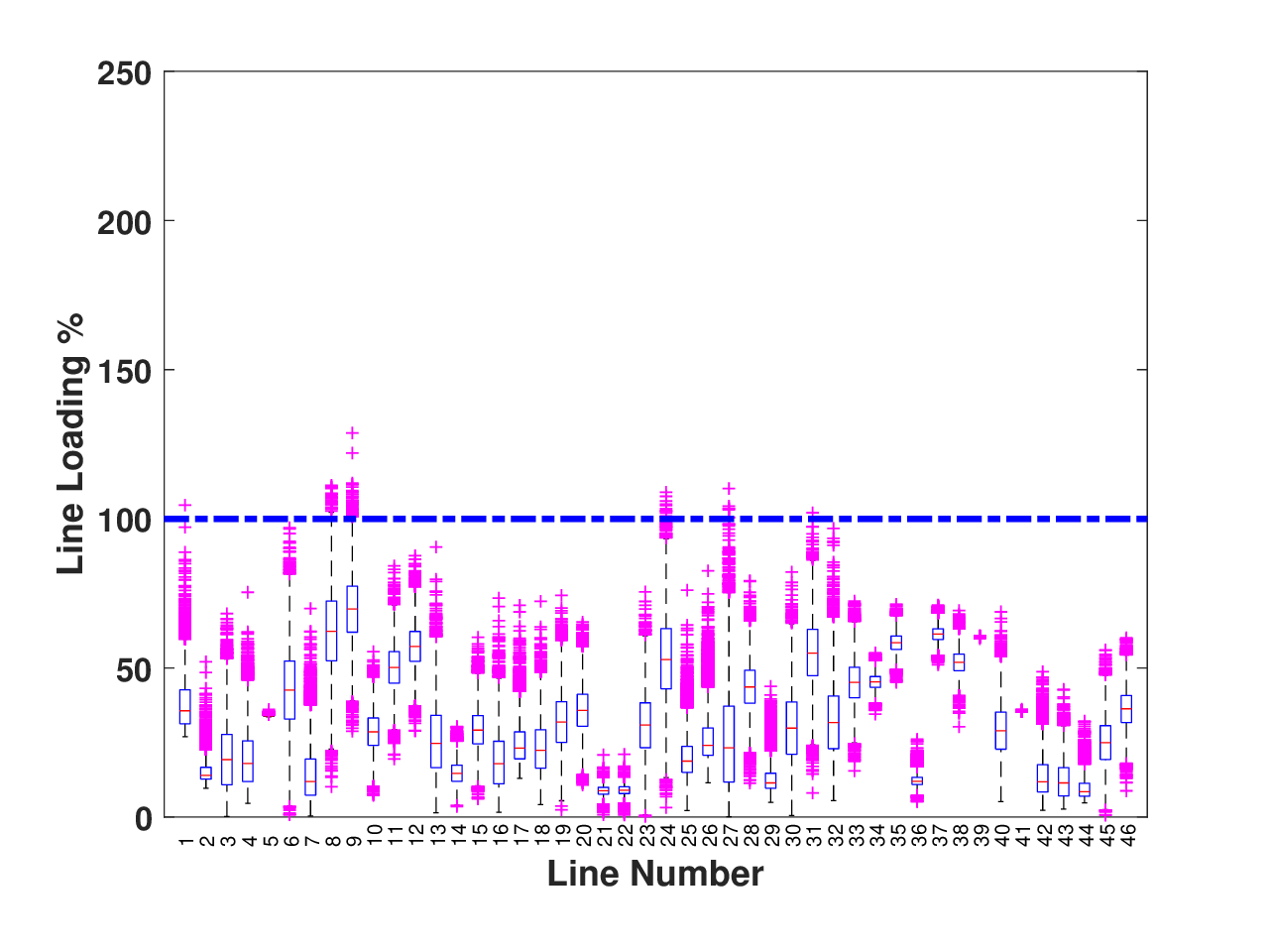}
		\caption{}
	\end{subfigure}
	\begin{subfigure}{.5\textwidth}
		\centering
		\includegraphics[trim={.3cm 1.1cm 1.2cm 1cm},clip,width=0.9\columnwidth]{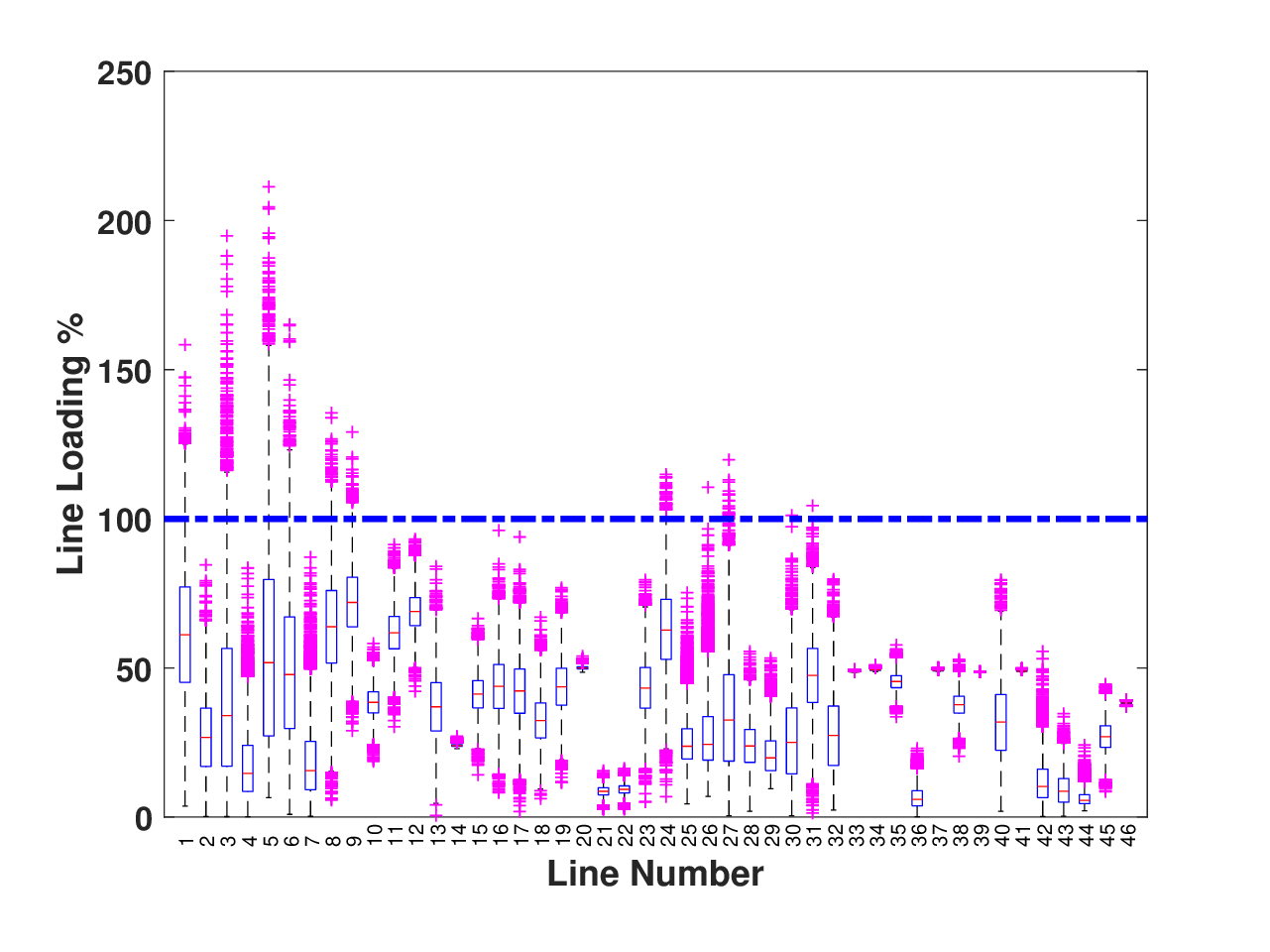}
		\caption{}
	\end{subfigure}
	\caption{\footnotesize (a): Line loading percentage for $10,000$ random samples extracted from the uncertainty set, results of scenario with certificates approach, (b) results of standard OPF design where no uncertainty is taken into account in designing the control parameters; 
    In the boxplots, the red line represents the median value,  edges of each box correspond to the $25^\text{th}$ and $75^\text{th}$ percentiles, whereas the whiskers extend to $99\%$ coverage. The magenta marks denote the data outliers. }
	\label{fig: OPF SImulation}
\end{figure*}

In Fig. \ref{fig: OPF SImulation}, we compare line loading percentage  for scenario with certificates approach against standard optimal power flow design where no uncertainty is taken into account while designing the control parameters. Demand level 
is $6,406$ MVA. This proves the necessity of adopting a robust strategy for the new power grids containing renewable  generators. In such a network if we rely on the classical OPF design where no uncertainty is taken into account in designing the control variables, it results in very frequent overloads, as shown in Fig. \ref{fig: OPF SImulation}(b), and hence frequent line tripping or even cascading outage. On the other hand, the robust strategy presented in the current paper successfully designs the control variables such that only very occasional violation happens during the operation of the network in the presence of large number of uncertain parameters.

We also compare the proposed strategy against the one presented in  \cite{ETH-AC}. In this comparison, we only consider uncertainty in renewable energy generators. Loads are considered to be known exactly. This is because  \cite{ETH-AC} can only handle uncertainty in renewable energy generator.
In \cite{ETH-AC}, the number of design variables appearing in the optimization problem is much larger than SwC approach. This is due to the linear parametrization  \eqref{eq:W-param} where some additional design variables are introduced. Since the number of scenario samples $N_\ped{SwC}$ depends on the number of design variables $n_\theta$, additional design variables lead to significant growth in the number of scenario samples as compared to SwC approach.
The worst-case cost is minimized in the OPF formulation---see subsection \ref{sec: worst case cost}. The penetration level is assumed to be $33\%$ and demand level is $7,112$ MVA corresponding to the peak demand in Fig. \ref{fig: 24hour results}.

Choosing $\varepsilon=0.1,~\beta=1\times10^{-10}$, the number of scenario samples for \cite{ETH-AC} is $11,838$---with $896$ design variables---and for scenario with certificates is $839$---with only $31$ design variables. Large number of scenario samples associated with  the approach presented in \cite{ETH-AC} results in much higher computational cost compared to SwC approach. 
Indeed, the optimization problem formulated based on \cite{ETH-AC} takes more than $51$ hours to be solved on a workstation with $12$ cores and $48$ GB of RAM, while the solution of the SwC approach  takes less than $49$ minutes. We remark that the computational time refers to a non-optimized implementation of the problem. 
 We also remark that the alternative approach proposed in \cite{ETH-AC}, and followed by recent works as \cite{Rostampour2017},
 based on computing a probabilistically guaranteed hyper-rectangle, may not be practically viable. Indeed, 
these require to solve a problem involving all the vertices of the uncertainty sets, whose number grows exponentially in the size of the uncertainty. For instance, in our numerical example,  it would amount at imposing the constraints on  $2^{46}\approx 7\times10^{13}$ vertices.
Finally, we note that sequential approaches as those proposed in \cite{Chamanbaz_TAC_2016}, and the 
sophisticated techniques of Sparsity Pattern Decomposition presented in  \cite{Rostampour2017},
may be adopted to significantly reduce computation times. 

\begin{figure}
\centering
\includegraphics[width=.6\columnwidth]{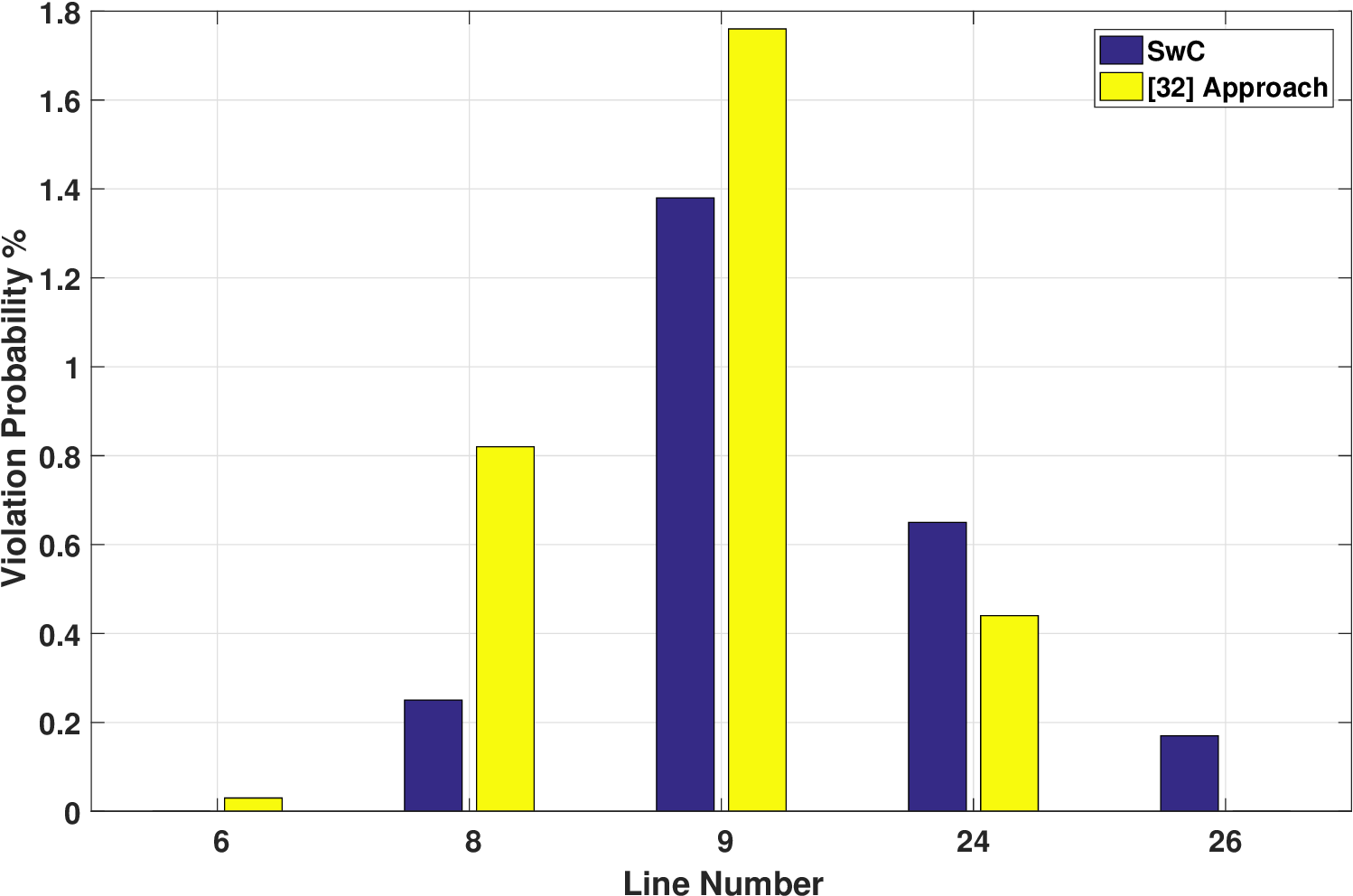}
\caption{\footnotesize The probability of line flow violation compared with \cite{ETH-AC}. The number of scenario samples for SwC is $839$ while \cite{ETH-AC} requires $11,838$ samples. }
\label{fig: ETHCompare}
\end{figure}

\begin{table}[!t]
\caption{\footnotesize Comparison between the average generation cost and average computation time -- over $10,000$ Monte Carlo simulation -- achieved  using SwC with the one obtained by~\cite{ETH-AC}.}
\begin{center}
\begin{tabular}{c||c||c}
\toprule
Approach &Generation Cost [\$]& Computation Time [min]\tabularnewline
\midrule
SwC & $25,280$ & 49\tabularnewline
\midrule
\cite{ETH-AC} & $25,359$ & 3,074\tabularnewline
\bottomrule
\end{tabular}
\end{center}
\label{tab: comarison of average cost}
\end{table}
Finally, we run an \textit{a-posteriori} analysis based on Monte Carlo simulation to estimate the probability of line flow violation.  In the posteriori analysis,  we use exactly the same set of samples---different from the design samples of course---to evaluate performance of the two approaches.  Figure \ref{fig: ETHCompare} 
compares the probability of line flow violation for the two approaches. 
This comparison shows that the linear parameterization adopted in \cite{ETH-AC} is conservative leading to larger violation level compared to SwC approach. The average---over $10,000$ samples---generation cost is also compared in Table \ref{tab: comarison of average cost}. Therefore, the paradigm presented in this paper improves upon \cite{ETH-AC} in computational complexity, violation probability, and generation cost.

\section{Concluding Remarks}\label{sec: conclusions}
In this paper, we proposed a novel approach to the AC optimal power flow problem in the presence of uncertain renewable energy sources and uncertain load. Assuming that the probability distribution of the uncertainty is available, we aim at optimizing the nominal power generation  subject to a small well defined risk of violating generation and transmission constraints. To tackle this complex NP-hard problem, we propose a randomized algorithm based on the novel concept of scenario with certificates and on convex relaxations of power flow problems. The effectiveness of the proposed solution is illustrated via numerical examples, where it is shown that one can significantly decrease the probability of constraint violation without a significant impact on the nominal power generation cost. Moreover, the approach is shown to be very efficient from a computational viewpoint. 
Future research aims at exploiting the flexibility of the methodology to extend its application.  
For instance, we can consider the case of optimizing with respect to the slack bus voltage, to improve the objective function and even the safety. Since the slack bus voltage is a discrete variable, the resulting OPF problem will be of mixed-integer type, as in \cite{MoGaLi:17}.

\bibliographystyle{plain}

\appendices
\section{Proof of Theorem~\ref{theo-OPF-robust}}\label{sec: appendix}
\noindent
Let $P_{k,opt}^G$ $\alpha_{opt}$ and $\mathbf{W}_{\mathbf{u},opt}$ be achievers of the solution of  {\bf CR-AC-OPF}. Now, take any $\dd\in\DD$ and consider the following optimization problem
{\small \begin{align*}
&\max_{\mathbf{W}\succeq 0} \  \sum_{k\in \mathcal{G}} Q_k \\
&\text{such that: } 
  P^G_{k,opt}+\alpha_{k,opt} \mathbf{s}^T \real\{\dd\}\leq  - P^R_{k}(\dd) +  P^L_{k}(\dd) +\\ &\sum_{l\in\mathcal{N}_k} \real\left\{(W_{kk}-W_{kl})^*y_{kl}^*\right\}   \leq
 P^G_{k,opt}+\alpha_{k,opt} \mathbf{s}^T \real\{\dd\},\,\,   \forall k \in \mathcal{N}\\
& Q_{k\,\min} \leq - Q^R_{k}(\dd) + Q^L_{k}(\dd) =\\ &\qquad\sum_{l\in\mathcal{N}_k} \imag\left\{(W_{kk}-W_{kl})^*y_{kl}^*\right\} \leq Q_{k\,\max}, \quad \forall k \in \mathcal{N}\\
& (\mathbf{W}^{\mathbf{u}}_{opt})_{kk} \leq W_{kk} \leq (\mathbf{W}^{\mathbf{u}}_{opt})_{kk} , \quad \forall k \in \mathcal{G};\\
& \left(V_{k\,\min}\right)^2\leq W_{kk}\leq \left(V_{k\,\max}\right)^2,  \quad \forall k \in \mathcal{N}/\mathcal{G}\\
& 
W_{ll}+W_{mm}-W_{lm}-W_{ml}\le (\Delta V_{lm}^{\max})^2,\ \quad\forall(l,m)\in\mathcal{L}.
\end{align*} } \vskip -.2in
Note that, for the given value of the uncertainty, the solution of the optimization problem above has exactly the same generation cost as  the solution of the {\bf CR-AC-OPF}. 
The optimization problem above is of the same form as that used in the proof of part (a) of Theorem~2 in~\cite{madani2015convex}. Hence, the same reasoning can be applied to show that there exists a rank one solution for the problem above. 
Since, as mentioned before, the solution of the convex relaxation and the one of the original optimal power flow problem coincide when $\text{rank}(\WW)=1$, this implies that the values of control variables obtained using  {\bf CR-AC-OPF} are a feasible power allocation for any value of the uncertainty   $\dd\in\DD$. Moreover, for all $\dd\in\DD$, the generation cost is equal to the optimal value of  {\bf CR-AC-OPF.} 
This, together with the fact that the solution of {\bf CR-AC-OPF} is optimal for at least one $\dd^*\in\DD$, leads to the conclusion that  {\bf CR-AC-OPF} does provide the best worst-case solution or, in other words, the robustly optimal power generation allocation.
\hfill $\square$\\
\vskip -.15in

\vspace{-.5in}
\begin{IEEEbiography}
	[{\includegraphics[width=1in,height=1.25in,clip,keepaspectratio]{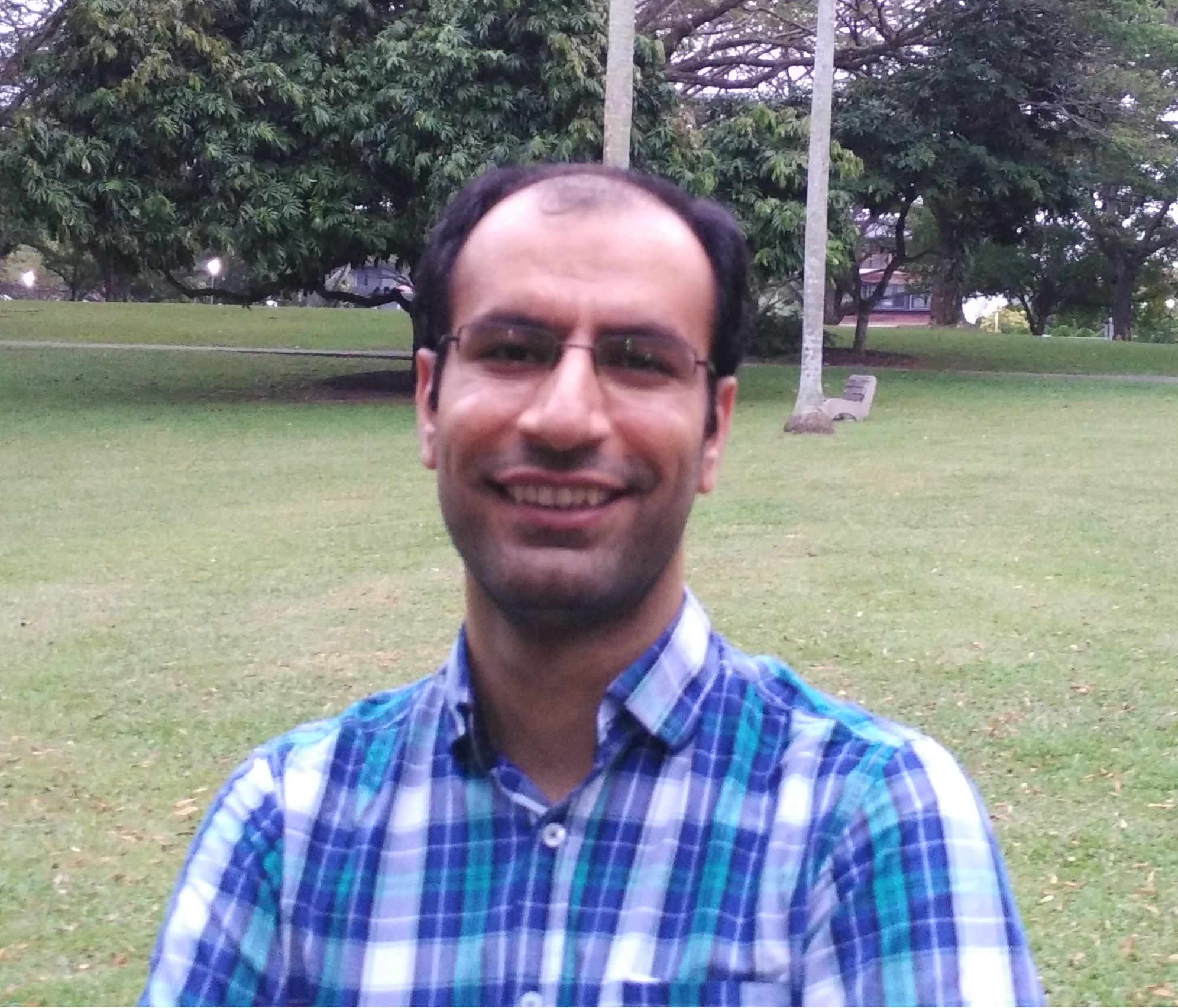}}]{Mohammadreza~Chamanbaz}  received his BSc in Electrical Engineering from Shiraz University of Technology. In 2014 he received his PhD from the Department of Electrical \& Computer Engineering, National University of Singapore in control science. Dr. Chamanbaz was with Data Storage Institute, Singapore as research scholar from 2010 to 2014. From 2014 to 2017, he was postdoctoral research fellow in Singapore University of Technology and Design. He was Assistant Professor in Arak University of Technology,  from Jan 2017 to Jan 2019. Dr. Chamanbaz is now senior research fellow in iTrust Center for Research in Cyber Security, Singapore.

His research activities are mainly focused on probabilistic and randomized algorithms for analysis and control of uncertain systems, robust and distributed optimization, and secure control of cyber-physical systems.

\end{IEEEbiography}
\vspace{-0.5in}	

\begin{IEEEbiography}
	[{\includegraphics[width=1in,height=1.25in,clip,keepaspectratio]{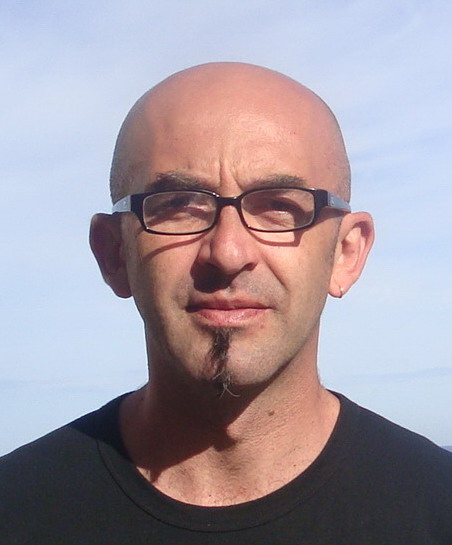}}]{Fabrizio~Dabbene} received the Laurea degree in 1995 and the Ph.D. degree in 1999, both from Politecnico di Torino, Italy. He is currently Senior Researcher at the CNR-IEIIT institute. His research interests include randomized and robust methods for systems and control, and modeling of environmental systems. He published more than 100 research papers and two books, and is recipient of the 2010 EurAgeng Outstanding Paper Award. He served as Associate Editor for Automatica (2008-2014) and IEEE Transactions on Automatic Control (2008-2012). Dr. Dabbene is a Senior Member of the IEEE, and has taken various responsibilities within the IEEE-CSS: he served as elected member of the Board of Governors (2014-2016) and as Vice President for Publications (2015-2016).
	
\end{IEEEbiography}
\vspace{-0.5in}

\begin{IEEEbiography}
	[{\includegraphics[width=1in,height=1.25in,clip,keepaspectratio]{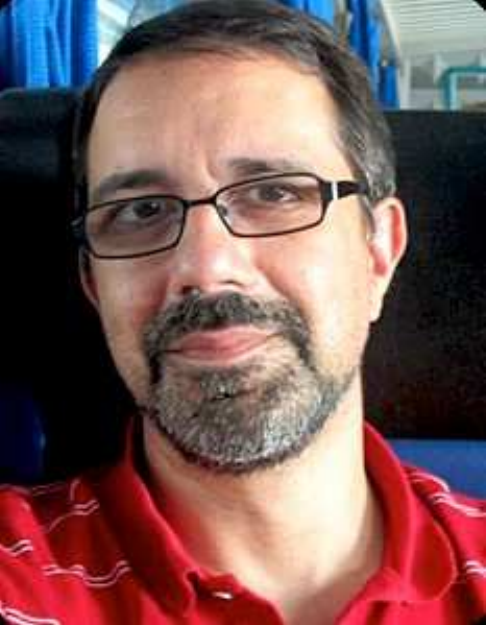}}]{Constantino~M.~Lagoa}
	received the B.S. and M.S. degrees from the
Instituto Superior Tecnico, Tech- ical University of Lisbon, Portugal
in 1991 and 1994, respectively, and the Ph.D. degree from the University
of Wisconsin at Madison in 1998. He joined the Electrical Engineering
Department of Pennsylvania State University, University Park, PA, in
August 1998, where he currently holds the position of Professor. He
has a wide range of research interests including robust optimization and
control, chance constrained optimization, controller design under risk
specifications, system identification and control of computer networks. Dr. Lagoa has served as Associate
Editor of IEEE Transactions on Automatic Control (2012-2017) and IEEE Transactions on Control systems Technology (2009-2013) and he is currently Associate Editor of Automatica.
	
\end{IEEEbiography}

\vfill

\end{document}